\documentclass[12pt]{article}

\usepackage{amsmath}
\usepackage{amssymb}

\usepackage[T1]{fontenc} 
\usepackage{sets pace}
\usepackage{comment}
\usepackage{colortbl}
\definecolor{text1}{cmyk}{1,.65,0,0} 
\definecolor{text2}{rgb}{1,0,0} 
\definecolor{text3}{cmyk}{0,0,0,1} 
\definecolor{text4}{cmyk}{0,0,0,0.5} 
\definecolor{text5}{cmyk}{1.0,0.0,1.0,0} 

\usepackage{amsfonts}
\usepackage{amsthm}
\usepackage{graphicx}
\usepackage{caption}
\usepackage{subcaption}
\usepackage{latexsym}
\usepackage[most]{tcolorbox}
\makeatletter
\renewcommand{\@seccntformat}[1]
{\csname the#1\endcsname.\enspace}
\makeatother
\newtheorem{theorem}{Theorem}
\newtheorem{lemma}{Lemma}
\newtheorem{remark}{Remark}

\newtheorem{example}{Example}

\begin{document}
\begin{center}
   {\bf Predictive density estimators with integrated $L_1$ loss \footnote{\today} }
\end{center}

\begin{center}
{\sc Pankaj Bhagwat$^{a}$ \& \'{E}ric Marchand
} \\

{\it a  Universit\'e de Sherbrooke, D\'epartement de math\'ematiques, Sherbrooke Qc,    CANADA, J1K 2R1 \quad 
(e-mails: pankaj.uttam.bhagwat@usherbrooke.ca; eric.marchand@usherbrooke.ca)} \\

\end{center}

\vspace*{0.2cm}
\begin{center}
{\sc Abstract} \\
\end{center}
\vspace*{0.1cm}
\small  
This paper addresses the problem of an efficient predictive density estimation for the density $q(\|y-\theta\|^2)$ of $Y$ based on $X \sim p(\|x-\theta\|^2)$ for $y, x, \theta \in \mathbb{R}^d$.   The chosen criteria are integrated 
$L_1$ loss given by $L(\theta, \hat{q}) \, =\, \int_{\mathbb{R}^d} \big|\hat{q}(y)- q(\|y-\theta\|^2)  \big| \, dy$, and the associated frequentist risk, for $\theta \in \Theta$.  For absolutely continuous and strictly decreasing $q$, we establish the inevitability of scale expansion improvements $\hat{q}_c(y;X)\,=\, \frac{1}{c^d} q\big(\|y-X\|^2/c^2  \big) $ over the plug-in density $\hat{q}_1$, for a subset of values $c \in (1,c_0)$.   The finding is universal with respect to $p,q$, and $d \geq 2$, and extended to loss functions $\gamma \big(L(\theta, \hat{q} ) \big)$ with strictly increasing $\gamma$.  The finding is also extended to include scale expansion improvements of more general plug-in densities $q(\|y-\hat{\theta}(X)\|^2 \big)$, when the parameter space $\Theta$ is a compact subset of $\mathbb{R}^d$.  Numerical analyses illustrative of the dominance findings are presented and commented upon.
As a complement, we demonstrate that the unimodal assumption on $q$ is necessary with a detailed analysis of cases where the distribution of $Y|\theta$ is uniformly distributed on a ball centered about $\theta$.  In such cases, we provide a univariate ($d=1$) example where the best equivariant estimator is a plug-in estimator, and we obtain 
cases (for $d=1,3$) where the plug-in density $\hat{q}_1$ is optimal among all $\hat{q}_c$.

\vspace*{0.2cm}


\normalsize
\section{Introduction}

We consider the problem of obtaining an efficient predictive density estimator $\hat{q}(y;X)$, $y \in \mathbb{R}^d$, of the density $q(\|y-\theta\|^2)$ of $Y$ based on spherically symmetric distributed $X \sim p(\|x-\theta\|^2)$.  In this set-up, the densities are Lebesgue on $\mathbb{R}^d$, $p$ and $q$ are known but not necessarily equal, and $X|\theta$ and $Y|\theta$ are independently distributed.  The observable $X$ may be a summary statistic arising from a sample.    We evaluate the efficiency of the predictive density of $\hat{q}(y;X)$ with integrated $L_1$ loss and risk
\begin{eqnarray}
\label{lossL1}
L(\theta, \hat{q}) \, & = &\, \int_{\mathbb{R}^d} \big|q(\|y-\theta\|^2) - \hat{q}(y)\big| \, dy \; , \\
R(\theta, \hat{q}) \, & = & \,  \int_{\mathbb{R}^d}   L\left(\theta, \hat{q}(\cdot; X)\right) \,\, p(\|x-\theta\|^2) \, dx\,.
\end{eqnarray} 
Spherically symmetric models are prominent in statistical theory and practice and inference for such models have a long history, including shrinkage estimation techniques (e.g., Fourdrinier et al., 2018).   Our set-up includes the normal case 
with 
\begin{equation}
\label{normalmodel}
X|\theta\sim N_d(\theta, \sigma^2_X I_d) \hbox{ and } Y|\theta \sim N_d(\theta, \sigma^2_Y I_d)\,,
\end{equation}
as well as scale mixtures of normals with $\sigma^2_X$ and $\sigma^2_Y$ random, and including multivariate Cauchy, Student, Laplace, Logistic distributions, and many others.
\begin{remark}
\label{remarkTV+overlap}
$L_1$ loss is a natural, appealing and widely used choice, which is also related to both:
\begin{enumerate}
\item[ (i)]
the ubiquitous total variation distance through the identity  $TV(f_1,f_2) \, = \, \sup_A \{\mathbb{P}(W_1 \in A) - \mathbb{P}(W_2 \in A)\} \, = \, \frac{1}{2} \, \int_{\mathbb{R}^d} | f_1(t) - f_2(t)| \, dt$, for random variables $W_1, W_2 \in \mathbb{R}^d$ with densities $f_1, f_2$; 
\item[ (ii)] the overlap coefficient (e.g., Weitzman, 1970) measuring the proximity of densities $f_1$ and $f_2$, and given by
\begin{equation}
\label{overlap}   OVL(f_1,f_2) \, = \, \int_{\mathbb{R}^d} \min(f_1(t), f_2(t)) \, dt \, = \, 1 - \frac{1}{2} \int_{\mathbb{R}^d} |f_1(t)-f_2(t)| dt \,,
\end{equation}
since $2\min(f_1, f_2) \, = \,f_1 + f_2 - |f_1-f_2| $.
\end{enumerate}
\end{remark}
 
There has been much interest over recent years and findings related to the efficiency of predictive density estimators in a decision-theoretic framework 
(e.g., George et al., 2019) and, in particular, relationships with shrinkage estimation techniques and Bayesian predictive densities.    However, frequentist $L_1$ risk predictive analysis is most challenging and the determination of Bayesian densities associated with loss (\ref{lossL1}), including the MRE predictive density obtained as a Bayes density with respect to the uniform measure $\pi(\theta)=1$, remains for the most part elusive  (but see Section \ref{sectionbayes} for exceptions).  Recently, Nogales (2021) considered loss $\{L(\theta, \hat{q})\}^2$ and showed that the Bayes predictive density matches quite generally the posterior predictive density (i.e., the conditional density of $Y$ given $X$ obtained by integrating out $\theta$).  It is thus of interest to study such loss functions as well, and more generally losses of the form $\gamma(L(\theta,\hat{q}))$.

Kubokawa et al. (2017) considered the benchmark choice for unimodal $q(\|y-\theta\|^2)$ of the plug-in predictive mle given by density $q(\|y-X\|^2)$.   Using an equivalence with a point estimation problem under a loss function which is a concave function of squared error loss, and shrinkage estimation techniques for such losses, they provided for $d \geq 4$, general $p$, and unimodal $q$, plug-in predictive densities of the form $q(\|y-\hat{\theta}(X)\|^2)$; $y \in \mathbb{R}^d$; that dominate the predictive mle under $L_1$ loss.    Their findings are quite general with respect to $p$ and $q$, but do require four dimensions or more.   Applications for normal and scale mixtures of normal distributions are expanded on as well. 

In the same paper, the authors provide and illustrate, for the univariate case ($d=1$), scale expansion improvements under $L_1$ loss of the type $\frac{1}{c} q(\frac{|y-x|^2}{c^2})$ for $c \in (1,c_0)$ on the plug-in predictive mle (i.e., $c=1$), requiring log-concavity of $q$.  For Kullback-Leibler loss, the potential of scale expansion improvements and potential 
inefficiency of a plug-in predictive density (also referred to an estimative fit) can be traced back to the work of Aitchison (1975), and is well illustrated by the 
normal case (\ref{normalmodel}) where the plug-in density $ \hat{q}_1 \sim N_d(X, \sigma^2_Y I_d)$, as well as all equivariant densities with respect to changes in location, are dominated by the MRE density $\hat{q}_U \sim N_d(X, (\sigma^2_X + \sigma^2_Y) I_d)$ which expands the scale indeed.  As a slightly tangential remark, but central to predictive density estimation findings over the past twenty years or so, we point out that $\hat{q}_U $ is under Kullback-Leibler loss minimax for all $d$, admissible for $d=1,2$, but inadmissible for $d \geq 3$ and dominated by various Bayesian predictive densities with striking parallels with shrinkage estimation under squared error loss and normal observables (Komaki, 2001; George et al., 2006; Brown et al. 2008).  

Returning to improvements by scale expansion, Fourdrinier et al. (2011) showed that any plug-in density $q(\|y-\hat{\theta}(X)\|^2)$ is dominated by a class of scale expansion variants for normal model $q$ and Kullback-Leibler loss.  Similar findings  were obtained by L'Moudden \& Marchand (2019) for normal models and $\alpha$-divergence loss, as well as by Kubokawa et al. (2015) for integrated $L_2$ loss.  Given the $L_1$ results for $d=1$ of Kubokawa et al. (2017) and their plug-in improvements on the predictive mle for $d \geq 4$, the remains the open questions of: {\bf (i)} improvements for $d=2,3$ and $\hat{\theta}(X)=X$ and {\bf (ii)} scale expansion improvements for $d \geq 2$ and $\hat{\theta}(X)=X$. 
We provide affirmative answers to these questions {\bf (i)} and {\bf (ii)}, as well as to: {\bf (iii)} scale expansion improvements on plug-in densities of the form $q(\|y-\hat{\theta}(X)\|^2)$ for choices of $\hat{\theta}(X)$ and when $\theta \in \Theta$ with $\Theta$ being compact, and even when $\hat{\theta}(X)$ is adapted to the parameter space.  Moreover, these results are established for the wider class of loss functions of the form $\gamma(L(\theta,\hat{q}))$ with strictly increasing $\gamma$.

\noindent The paper is organized as follows.  In {\bf Section 2}, 
we first expand on Bayesian predictive densities under $L_1$ loss and the general difficulty in determining a Bayesian solution, nevertheless recording an explicit solution for the univariate uniform distribution case.  {\bf Section 3} contains the main dominance findings which relates to predictive densities of the form 
\begin{equation}
\label{qhatc}
q_{\hat{\theta},c}(y;X) \, = \, \frac{1}{c^d} \, q(\frac{\|y-\hat{\theta}(X)\|^2}{c^2})\,, y \in \mathbb{R}^d \,.
\end{equation}
In {\bf Section 3.1}, we study cases with $\hat{\theta}(X)=X$, general $p$ and $q$ with unimodal $q$, and losses $\gamma(L(\theta,\hat{q}))$ with strictly increasing $\gamma$.  With such densities having  constant risk, we show that the optimal scale expansion value $c^*$ is such that $c^* >1$.   The proof is unified and applicable for quite generally for arbitrary $(p,q,d, \gamma)$ such that $d \geq 2$ and $q$ is strictly decreasing on $\mathbb{R}_+$.   Secondly in {\bf Section 3.2}, we consider situations with a compact parameter space restriction $\theta \in C$, such as balls of a  radius $m$, and show quite generally that a plug-in density $q_{\hat{\theta},1}(\cdot;X)$ is necessarily dominated by a subclass of scale expansion variants
$q_{\hat{\theta},c}(\cdot;X)$ with $c \in (1,c_0)$.   The finding is again unified for general $p$, unimodal $q$, $\gamma$, and $d \geq 2$.   We do also provide in {\bf Section 3.3} cases where the plug-in density is optimal among $q_{\hat{\theta},c}$, namely for uniformly distributed $Y \sim U(\theta-B, \theta+B)$.  We further explore such phenomena in the multivariate case with $Y$ uniformly distributed on a ball of radius $m$, centered at $\theta$ aided by numerical evaluations and a definite result for $d=3$ and $X$ uniformly distributed on a ball centered at $\theta$.  Finally, numerical illustrations and comparisons, as well as concluding observations, are presented in {\bf Sections 4 and 5}. 

\section{Bayesian predictive densities and $L_1$ loss}
\label{sectionbayes}

Despite the appeal of $L_1$ divergence for reporting on the efficiency of estimated densities in parametric and non-parametric settings (e.g., DasGupta \& Lahiri, 2012; Devroye \& Györfi, 1985), drawbacks include the challenging frequentist risk analysis and, mostly, the difficulty of specifying a Bayesian predictive density.   This difficulty includes the determination of the minimum risk equivariant density for location models or, equivalently, the Bayes predictive density $\hat{q}_{mre}$ with respect to the uniform density $\pi(\theta)=1$.  The equivalence follows from a general representation for the minimum risk equivariant estimator as the Bayes estimator associated with the corresponding Haar measure (e.g., Eaton, 1989) which is the uniform density derived from the group of location changes.   Moreover, given general results on equivariant procedures (Kiefer, 1957), such a density is minimax and thus constitutes an interesting benchmark predictive density.   

In this section, we briefly expand on such a difficulty and it is particularly instructive to contrast Bayesian solutions for Kullback-Leibler and $L_1$ loss functions.    Moreover, we do provide a Bayes predictive density solution under $L_1$ loss in the uniform case with unknown location (Theorem \ref{theoremuniform}).
Interestingly, defining the loss as the square of the $L_1$ distance leads to tractable Bayesian solutions (Nogales, 2021).

For the sake of illustration, we suppose in this section that $X \sim p_{\theta}$ is to be observed and that we wish to obtain a predictive density $\hat{q}(y ;X)$, $y \in \mathbb{R}^d$, for the density $q_{\theta}$ of $Y$.  We assume that  $p_{\theta}$ and  $q_{\theta}$ are Lebesgue densities, and that $X|\theta$ and $Y|\theta$ are independently distributed.   Finally, we assume a prior density $\pi$ for $\theta$ for which the posterior $\pi(\cdot|x)$; defined with respect to $\sigma-$finite measure $\nu$; exists.  
\subsection{Kullback-Leibler loss}

The familiar Kullback-Leibler loss associated with density $\hat{q}$ as an estimate of $q_{\theta}$ is given by
\begin{equation}
\nonumber
L_{KL}(\theta, \hat{q}) \, = \, 
\int_{\mathbb{R}^d}  q_{\theta}(y)  \log \{ \frac{q_{\theta}(y)}{\hat{q}(y)} \} \, dy\,.
\end{equation}
A useful and equivalent representation is given by 
\begin{equation}
\label{lossKLalternative}
L_{KL}(\theta, \hat{q}) \, = \, 
\int_{\mathbb{R}^d}  q_{\theta}(y)  \left\lbrace \frac{\hat{q}(y)}{q_{\theta}(y)} - \log ( \frac{\hat{q}(y)}{q_{\theta}(y)} ) -1 \right\rbrace \, dy \,.
\end{equation}

The above clearly represents the loss as a weighted (with respect to $q_{\theta}(y)$) average of a collection of
distances between estimates $\hat{q}(y)$ and actual $q_{\theta}(y)$ values as measured by the point estimation loss $\rho(\frac{\hat{q}(y)}{q_{\theta}(y)})$ with $\rho(z)=z-\log z -1$.   

Now, consider estimating the density $q_{\theta}(t)$ at a fixed value $y=t$ and refer this as the local problem.    The Bayes estimate minimizes in $\hat{q}(t)$ the expected posterior loss

$$\mathbb{E} \{ q_{\theta}(t) \, \rho(\frac{\hat{q}(t)}{q_{\theta}(t)})|x\}$$
 in $\hat{q}(t)$.  It is then easy to infer that the local Bayes estimate is given by 
\begin{equation}
\label{bayeskl} 
 \hat{q}_{\pi}(t;x) \, = \, \mathbb{E} (q_{\theta}(t)  |x) \, = \, \int_{\Theta} q_{\theta}(t) \, \pi(\theta|x) \, d\nu(\theta)\,. 
 \end{equation}

Now, for the global problem, a Bayesian predictive density $\hat{q}(y;x)$, $y \in \mathbb{R}^d$, minimizes among all densities the expected posterior loss which, from (\ref{lossKLalternative}) and a change in the order of integration, becomes equivalent to minimizing 
\begin{equation}
\int_{\mathbb{R}^d}    \mathbb{E} \{ q_{\theta}(y) \, \rho(\frac{\hat{q}(y)}{q_{\theta}(y)})|x\} \; dy\,.
\end{equation}
Finally, since $\hat{q}_{\pi}(y;x)$ minimizes for all $y$ the expectation inside the above integral, and since $\hat{q}_{\pi}(\cdot;x)$ is actually a density on $\mathbb{R}^d$, it follows that $\hat{q}_{\pi}(\cdot;x)$ is the Bayes predictive density. 

\subsection{$L_1$ loss}

The approach presented for KL loss is possible basically since a probabilistic weighted average of densities is a density.    For integrated $L_1$ loss, the local problem brings into play the median value of $q_{\theta}(t)$ with respect to the posterior distribution of $\theta$.  However in general, such a collection of median values do not form a density and the global minimization problem cannot be deduced from the local problems.   Nevertheless, we next record such a possibility where the resulting values do form a density and follow-up with an application to the uniform model with unknown location parameter and an explicit expression for $\hat{q}_{mre}$.  To proceed, we denote
\begin{equation} 
\label{medianposterior}
 \hat{q}(t;x) \, =\, \hbox{Med}  \{q_{\theta}(t) |x\} \,, x, t \in \mathbb{R}^d\,,
\end{equation}
which is not necessarily unique.
\begin{lemma}
\label{lemmamedian}
If, for the general set-up of this section, $\hat{q}(\cdot;x)$ is a density on $\mathbb{R}^d$ for all $x$, then it is a Bayes predictive density with respect to $L_1$ loss.
\end{lemma}
{\bf Proof.}   This follows as in Section 2.1.  \qed

\begin{theorem}
\label{theoremuniform}
Let $X=(X_1, \ldots, X_n)$ with independently distributed $X_i \sim U(\theta-A, \theta+A)$, and independently of $Y \sim U(\theta-B, \theta+B)$.   Then, whenever $B \geq A/2$, a Bayesian predictive density associated with the uniform density $\pi(\theta)=1$ and under $L_1$ loss in (\ref{lossL1}) is given by a $U(\frac{X_{(1)}+X_{(n)}}{2} -B, \frac{X_{(1)}+X_{(n)}}{2} +B)$ density, with $X_{(1)}=\min \{ X_1, \ldots, X_n \}$ and $X_{(n)}=\max \{ X_1, \ldots, X_n \}$. 
\end{theorem}
{\bf Proof.}   Since $X_1,\ldots, X_n| \theta \hbox{ are i.i.d. } U(\theta - A,\theta + A)$ and $\pi(\theta) = 1$, we have the posterior distribution $\theta|x \sim U(x_{(n)} - A, x_{(1)} + A)$ with c.d.f.
\begin{equation}
\mathbb{P}(\theta \leq u | x)=\begin{cases}
          0 \quad &\text{if} \, u <  x_{(n)}  - A \\
\label{cdfc}          \frac{u-x_{(n)} + A}{x_{(1)} - x_{(n)} + 2A} \quad &\text{if} \, u \in [x_{(n)} - A,x_{(1)} + A) \\
          1 \quad &\text{if} \, u \geq  x_{(1)}  + A \,.\\
     \end{cases}
     \end{equation}

With the density $q_{\theta}(y)$ taking the values $\frac{1}{2B}$  and $0$ only, we have that for fixed $y$,
\begin{align}
  \nonumber    Median\{q_{\theta}(y)|x \} & = \frac{1}{2B} \\
\nonumber  \Longleftrightarrow  \quad  \mathbb{P}(q_{\theta}(y)= \frac{1}{2B}|\,x ) & \geq \frac{1}{2} \\
\nonumber \Longleftrightarrow \quad   \mathbb{P}(\theta \leq y + B|x) & - \mathbb{P}(\theta \leq y - B|x)  \geq \frac{1}{2}\,.
\end{align}

An evaluation of this posterior probability using (\ref{cdfc}) deploys itself into cases  : ({\bf I}) $B \geq A$ and ({\bf II}) $A > B \geq \frac{A}{2}$.
For ({\bf I}), we have
\begin{equation}
\mathbb{P}(\theta \leq y + B|x) - \mathbb{P}(\theta \leq y - B|x) =\begin{cases}
          0 \quad &\text{if} \, y <  x_{(n)}  - A  - B \text{ or } y > x_{(1)} + A + B \\
          \frac{y + B + A - x_{(n)} }{x_{(1)} - x_{(n)} + 2A} \quad &\text{if} \, y \in [x_{(n)} - A - B,x_{(1)} + A- B)\\
          1 &\text{if} \, y \in [x_{(1)} + A - B,x_ {(n)} - A + B] \\
         \nonumber 
                
          \frac{x_{(1)} + B + A - y }{x_{(1)} - x_{(n)} + 2A} \quad &\text{if} \,  y \in  (x_{(1)}  - A + B < y <   x_{(n)}  + A  + B )\, . \\
     \end{cases}
     \label{cdf}
\end{equation}

From this, we obtain that
\begin{align*}
    \mathbb{P}(\theta \leq y + B|x) - \mathbb{P}(\theta \leq y - B|x) \geq \frac{1}{2} \Longleftrightarrow y \in \left[ \frac{x_{(1)} + x_{(n)}}{2} - B,  \frac{x_{(1)} + x_{(n)}}{2} + B \right],
\end{align*}

which implies indeed that $Median\{q_{\theta}(y)|x \} \, = \, \frac{1}{2B} \,
\mathbb{I}_{\left(-B + \frac{x_{(1)} + x_{(n)}}{2} \,, \,  B+ \frac{x_{(1)} + x_{(n)}}{2}  \right)}(y)$, and that $\hat{q}_{\pi}(y;X)$ is a $U(\frac{X_{(1)}+X_{(n)}}{2} -B, \frac{X_{(1)}+X_{(n)}}{2} +B)$ density.  Finally, case ({\bf II}) is handled in a similar fashion, is left to the reader, and leads to the result.  \qed

We pursue with various observations.

\begin{remark}
As a complement to the above Theorem, consider the case $B < \frac{A}{2}$.  For such a case whenever $x_{(n)}-x_{(1)} < 2A - 4B$, a calculation shows that
\begin{align*}
    \max_{y \in \mathbb{R} } \mathbb{P}(\theta \in (y - B, y + B)|x) = \frac{2B}{x_{(1)} - x_{(n)} + 2A} \,,
\end{align*}
which implies that $Median\{q_{\theta}(y)|x\} = 0 $ for all $y$, and which obviously does not lead to valid density estimates for such $x$ values.  
\end{remark}

\begin{remark}
The above Bayesian predictive density is a plug-in density with $\hat{\theta}(X) \, = \, (X_{(1)}+X_{(n)})/2$ and the result also applies for $n=1$.  The estimator $\hat{\theta}(X)$ is familiar and plausible.  Namely, it is best unbiased  with respect to loss $\rho(\hat{\theta}-\theta)$ with $\rho$ even and convex as long as $\mathbb{E}(|\rho(\hat{\theta}(X)|) < \infty$.  From a Bayesian perspective, it matches the posterior mean and median associated with prior density $\pi(\theta)=1$.
\end{remark}

\begin{remark}
\label{remarkno}
We point out that Theorem \ref{theoremuniform}'s minimum risk equivariant predictive density for $Y$ has the same standard deviation as the model density for $Y,$ so that there is no beneficial scale expansion as such expansions are also equivariant, and equivariant decision rules here have constant risk.  The optimality of the plug-in mle among scale variants will be extended in Section \ref{uniformscaleexpansion} for other model distributions for $X$.
\end{remark}

\section{Risk analysis and dominance findings}

This section contains our main results concerning the effect of scale expansion in constructing improved density estimators for spherically symmetric models.  We consider $X \sim p(||x - \theta||^2)$ and $Y \sim q(||y - \theta||^2)$, independent conditional on $\theta \in \mathbb{R}^d$.   Section 3.1 relates to dominating the predictive mle for $d \geq 2$, while further dominance results applicable to plug-in densities, and to when $\theta$ is restricted to a compact subset of $\mathbb{R}^d$, are given in Section 3.2.  These findings relate to unimodal $q$, and Section 3.3's risk analysis for the case where the target density is that of a uniform distribution on balls around $\theta$ provide ``counterexamples'' where no beneficial scale expansion is possible.  Namely, for the uniform univariate case with $p$ non-increasing on $\mathbb{R}_+$, we show that the plug-in density is optimal among scale variants $q_{\hat{\theta},c}$ as defined in (\ref{qhatc}).

\subsection{Dominating the predictive mle density}

We assume here that $p$ and $q$ are Lebesgue densities and that $q$ is absolutely continuous and decreasing on $\mathbb{R}_+$.  Our target predictive density is the predictive mle which is simply the plug-in density $q(\|y-X\|^2)$, $y \in \mathbb{R}^d$.  We investigate the frequentist risk performance associated with 
$L_1$ loss of the following class of scale expansions densities $\hat{q}_{c}(y;X)  = \frac{1}{c^d} q\left(\frac{||y - X||^2}{c^2}\right)\,, y \in \mathbb{R}^d$.  The main finding is the inevitability of dominating $\hat{q}_c$'s for sufficiently small $c$, i.e., $1 < c  \leq c_0$ for some $c_0$ depending on $d, p, q$.   The result is also extended to loss functions $\gamma(L(\theta, \hat{q}))$ with strictly increasing $\gamma$.

We will make use of the following key result.   In the following, we denote   
$S_r$ as the $d-$dimensional unit sphere centered at $0$ given by $S_{r} = \left\{u \in \mathbb{R}^d \;:\; \|u\| = r\right\}$.  
\begin{lemma}
\label{lemmarandomunit}
Let $d \geq 2$, $U_1$ and $U_2$ be independent such that $U_1$ is uniformly distributed on $S_1$, and $U_2$ a random unit vector such that $\mathbb{P}(U_2 \in S_1)=1$.   Define $V = U_1^{\top}U_2$.    Then, the distribution of $V$ is independent of that of $U_2$ and has p.d.f.
\begin{equation}
\label{densityV}
f_{V}(v) = \, \frac{1}{B\left(\frac{d - 1}{2},\frac{1}{2}\right)} \,\, (1 - v^2)^{\frac{d - 3}{2}} \;,\; v \in [-1,1]\,,
\end{equation}
where ${\displaystyle B}(\cdot, \cdot)$ is the beta function.  
\end{lemma} 
{\bf Proof.}  Density (\ref{densityV}) is known to describe the distribution of $a^{\top}U_2$ for fixed $a \in S_1$ (e.g., Kariya \& Eaton, 1977).   The result therefore follows by conditioning and the independence assumption.  \qed 

\begin{theorem}
\label{theoremmain}
Let $X \sim p(||x - \theta||^2)$ and $Y \sim q(||x - \theta||^2)$, $x,y,\theta \in \mathbb{R}^d$, with $d \geq 2$, $q$ absolutely continuous and strictly decreasing on $\mathbb{R}_+$.   Consider estimating $q(||y - \theta||^2)$ based on $X$ under $\gamma(L(\theta, \hat{q}))$ loss, with $\gamma$ a strictly increasing, absolutely continuous function such that $\gamma(2)<\infty$, and predictive density estimators $\hat{q}_c(y;X) \, = \, \frac{1}{c^d}q\left(\frac{||y - X||^2}{c^2}\right)\,, y \in \mathbb{R}^d$.  Then
$\hat{q}_1(\cdot;X)$ is inadmissible and dominated by $\hat{q}_c(\cdot;X)$
for $c \in (1,c_0)$ and some $c_0> 1$, and of which there exists an optimal $\hat{q}_{c^*}$. 
\end{theorem}

{\bf Proof.}  The proof is divided into two parts {\bf (A)} and {\bf (B)}, first showing that $\hat{q}_c(\cdot;X)$ has constant risk $R(\theta, \hat{q}_c)$ in $\theta \in \mathbb{R}^d$ and given by:
\begin{equation}
\label{R(c)gamma}
R(c) = \, \mathbb{E}_0 \gamma \big[ 2 \mathbb{E}_0 \left\{F_V \big( l_1^c(||X||,||Y||)\big) - F_V \big( l_2^c(||X||,||Y||)\big)\right\}\big],
\end{equation}
where (i) $F_V$ is the c.d.f. associated to the p.d.f. in (\ref{densityV}), (ii) the outside and inside expectations $\mathbb{E}_0$ are taken with respect to the densities $p(\|x\|^2)$ and $q(\|y\|^2)$, respectively, and (iii) $l_1^c(t_1,t_2)$  and $ l_2^c(t_1,t_2) $ are to defined below.   Part {\bf (B)} will simply consist in showing that $R(c)$ decreases locally at $c=1$.

{\bf (A)}
The risk function under $L_1$ loss of $\hat{q}_{c}(\cdot;X)$ is given by
\begin{align}
\nonumber  R(\theta,\hat{q}_{c}) & = \int\limits_{\mathbb{R}^d} \gamma \big\{ \int\limits_{\mathbb{R}^d}\, \big|\,q(||y - \theta||^2) - \frac{1}{c^d} q\left(\frac{||y - x||^2}{c^2}\right)\,\big|\; dy \big\} \; p(||x - \theta||^2) \; dx \\
\label{formularisk}
& = \int\limits_{\mathbb{R}^d} \gamma \big\{\int\limits_{\mathbb{R}^d}\,\big|\,q(||y||^2) - \frac{1}{c^d} q\left(\frac{||y - x||^2}{c^2}\right) \big| \; dy \big\} \; p(||x||^2) \; dx \,
\\
\nonumber & = R(c) \,\; (\hbox{independently of } \theta),
\end{align}
with the change of variables $(x,y) \rightarrow (x - \theta, y -\theta)$.  Setting $A_1(x,c) =  \left\{y \in \mathbb{R}^d\;:\; q(||y||^2) \;  \geq \frac{1}{c^d} \; q\left(\frac{||y - x||^2}{c^2}\right)  \right\}$,
the above decomposes as 
\begin{eqnarray}
\nonumber R(c) & = & \int\limits_{\mathbb{R}^d} \gamma \bigg[\int\limits_{A_1(x,c)}
\, \big\{q(||y||^2) - \frac{1}{c^d} q(\frac{||y - x||^2}{c^2}) \big\} \;dy \,   \\
\nonumber \, & + & \int\limits_{\mathbb{R}^d} \, \big\{\frac{1}{c^d} q(\frac{||y -x||^2}{c^2})  \, - \, q(||y||^2) \big\} \;dy - \int\limits_{A_1(x,c)} \, 
\big\{\frac{1}{c^d} q\left(\frac{||y - x||^2}{c^2}\right)  \, - \, q(||y||^2) \big\} \;dy \bigg] \, p(||x||^2) \; dx  \\
\nonumber & = & \int\limits_{\mathbb{R}^d} \gamma \bigg\{ 2 \int\limits_{A_1(x,c)}\,q(||y||^2) \; dy \;  -  2 \int\limits_{A_1(x,c)}\, \frac{1}{c^d} q\left(\frac{||y - x||^2}{c^2}\,\right)\; dy \bigg\}\; p(||x||^2) \; dx \\
\nonumber & = & \int\limits_{\mathbb{R}^d} \gamma \bigg\{ 2 \int\limits_{A_1(x,c)}\,q(||y||^2) \; dy \;  -  2 \int\limits_{\mathcal{A}_2(x,c)}\,q(||y||^2)\; dy \bigg\}\; p(||x||^2) \; dx \, \\
 & = & \mathbb{E}_0 \, \gamma\{2\mathbb{P}_0 \big(Y \in A_1(X,c)\big) \, - \,  2 \, \mathbb{P}_0 \big(Y \in \mathcal{A}_2(X,c)\big)\} \,,
 \label{Rc2}
\end{eqnarray}
where $\mathcal{A}_2(x,c)$ is the image of $A_1(x,c)$ under the transformation $\frac{y - x}{c} \to y$ defined as $\mathcal{A}_2(x,c) =  \left\{y \in \mathbb{R}^d\;:\; q(||cy + x||^2) \;  \geq \frac{1}{c^d} \; q\left(||y||^2\right)  \right\}$,
the expectation is taken with respect to $X$ at $\theta=0$, and $\mathbb{P}_0$ is taken with respect to $Y$ (which is independent of $X$) at $\theta=0$. 
We now define the generalized inverse of $q$ as $q^{-1}(t)\,=\,\inf\{z \geq 0: q(z) \leq t\}$ in such a way so that $q^{-1}(t)=0$ whenever $t \geq q(0)$.
Since $q$ is monotonically decreasing on $\mathbb{R}_{+}$, we have 
\begin{align}
\nonumber A_1(x,c) =  & \left\{y \in  \mathbb{R}^d\;:\; q(||y||^2) \; c^d  \geq \; q\left(\frac{||y - x||^2}{c^2}\right)  \right\}   = \left\{y \in  \mathbb{R}^d \;:\; q^{-1}( q(||y||^2) \;c^d  ) \leq \; \frac{||y - x||^2}{c^2}  \right\} \\
\nonumber  =&   \left\{y \in  \mathbb{R}^d \;:\;  \frac{y^{\top}x}{||y||||x||} \leq l_{1}^c(||x||,||y||)  \right\}\,,
\end{align}
where  $l_1^c(t_1,t_2)  =  \frac{t_1^2 + t_2^2 - c^2 q^{-1} \big( q(t_2^2) \; c^d \big)}{2t_1t_2} \textit{ for } t_1,t_2 > 0.$  Similarly, we get 
\begin{equation}
\nonumber
\mathcal{A}_2(x,c)  =  \left\{y \in \mathbb{R}^d\;:\;  \frac{y^{\top}x}{||y||||x||} \leq l_{2}^c(||x||,||y||)  \right\} , 
\end{equation}
where $l_2^c(t_1,t_2)  =  \frac{-t_1^2 -c^2 t_2^2 +  q^{-1}\left( q(t_2^2)/c^d \right) }{2\, c\, t_1 \, t_2} \textit{ for } t_1,t_2 > 0.$

Now, using Lemma \ref{lemmarandomunit}, it follows that $V \equiv \frac{Y^{\top}X}{||Y||||X||}$ is distributed independently of $(\|X\|, \|Y\|)$ with c.d.f. $F_V$.  Therefore, for $i=1,2$, we have
\begin{eqnarray}
\nonumber
\mathbb{P}_0\left( Y \in A_i(X,c) \right|X) & = &
 \mathbb{P}_0\left(V \leq l_{i}^c(||X||,||Y||) \right) \\
& = & \mathbb{E}_0 \left( F_V(l_{i}^c(||X||,||Y||))  \right) \,,
\end{eqnarray}
and (\ref{R(c)gamma}) thus follows from the above and (\ref{Rc2}).

{\bf (B)} It suffices to show that $\frac{\partial}{\partial c}R(c) \Bigr|_{\substack{c \rightarrow 1^{+}}} < 0 $. Differentiating $R(c)$ in (\ref{R(c)gamma}) under the integral sign, we get
\begin{equation}
\label{R'c}  
\frac{\partial}{\partial c}R(c) = 2\mathbb{E}_0^X  \{\gamma' (H_c(\|X\|^2) \,\mathbb{E}_0^Y G_c(||X||,||Y||)     \}\,,
\end{equation}
with $H_c(\|X\|^2) \,=\, 2 E_0^Y \{F_{V}(l_1^c(||X||,||Y||) \, - \, F_{V}(l_2^c(||X||,||Y||)\}$ and $$G_c(||X||,||Y||) \, = \, \left\{f_{V}(l_1^c(||X||,||Y||))\frac{\partial}{\partial c}l_1^c(||X||,||Y||) - f_{V}(l_2^c(||X||,||Y||))\frac{\partial}{\partial c}l_2^c(||X||,||Y||)\right\}\,. $$

We have
\begin{align}
\nonumber \frac{\partial}{\partial c}\mathnormal{l}_1^{c}(t_1,t_2) & = \frac{-2c q^{-1}( q(t_2^2) \;c^d ) -  dc^{d+1} \frac{ q(t_2^2)}{q'\{q^{-1}( q(t_2^2) \;c^d )\}}}{2 t_1 t_2}\,, \\
\nonumber \frac{\partial}{\partial c}\mathnormal{l}_2^{c}(t_1,t_2) & = \frac{\frac{1}{c^2} t_1^2 - t_2^2    - \frac{1}{c^2}q^{-1}( \frac{q(t_2^2)}{c^d} ) -  \frac{d}{c^{d+2}} \frac{\, q(t_2^2)}{\, q'\{q^{-1}(\frac{q(t_2^2)}{c^d}) \}}}{2 t_1 t_2 }\,,
\end{align}
$l_1^1(t_1,t_2) = -l_2^1(t_1,t_2) = \frac{t_1}{2t_2}$, and therefore
\begin{align*}
\lim_{c \to 1^+}\frac{\partial}{\partial c}\mathnormal{l}_1^{c}(t_1, t_2) = \frac{-2 t_2^2 - \frac{dq(t_2^2)}{q'(t_2^2)}}{2 t_1 t_2}\; \text{ and }\; \lim_{c \to 1^+} \frac{\partial}{\partial c}\mathnormal{l}_2^{c}(t_1,t_2) =  \frac{t_1^2 - 2 t_2^2 - \frac{dq(t_2^2)}{q'(t_2^2)}}{2 t_1 t_2} \,.
\end{align*}
Finally, we obtain from (\ref{R'c})
\begin{align*}
\frac{\partial}{\partial c}R(c) \Bigr|_{\substack{c \rightarrow 1^{+}}} =  -\mathbb{E}_0^X  \big\{\gamma' (H_1(\|X\|^2)) \,  \mathbb{E}_0 f_V\left(\frac{||X||}{2||Y||}\right)\frac{||X||}{||Y||} \big\} <  0\,, 
\end{align*}establishing the result.  \qed

We point out that the above finding and proof are unified for all unimodal $q$, all $p$, dimension $d \geq 2$, and choice of $\gamma$.   In fact as seen with the proof, the result applies to all spherically symmetric distributed $X$ as long as there is no atom at the origin (i.e., $\mathbb{P}_0(X=0)=0$).  This extends Kubokawa et al.'s (2017) $d=1$, identity $\gamma$ scale expansion improvement to $d \geq 2$.   It also establishes the inadmissibility of $\hat{q}_1(\cdot;X)$ for $d=2,3$, which was established by Kubokawa et al. (2017) for $d \geq 4$.   Combined with previous work, improvements for identity $\gamma$ on the predictive mle density thus arise with either scale expansion or again plug-in improvements of the type $q(\|y-\hat{\theta}(X)\|^2)$, $y \in \mathbb{R}^d$.  A detailed illustration is presented in Section 4.

\subsection{Dominating a plug-in density}

We now study the frequentist risk performance of densities $q_{\hat{\theta},c}$ as defined in (\ref{qhatc}) for a given estimator $\hat{\theta}(X)$ of $\theta$ and varying $c$.  Section 3.1 sets $\hat{\theta}(X)=X$, but we consider here more general non-degenerate choices of $\hat{\theta}(X)$, namely in the context of a compact parameter space restriction, such as balls of a fixed radius.  The objective remains to assess whether or not the plug-in density $q_{\hat{\theta},1}$ is improvable by a scale expansion variant $q_{\hat{\theta},c}$ with $c>1$.   Whereas densities $\hat{q}_c$ have constant risk, facilitating the risk analysis,  this will not be the case for different choices of $\hat{\theta}(X)$.   However, the following adaptation of Theorem \ref{theoremmain} which capitalizes on the compactness of the parameter space leads to the following dominance result.

\begin{theorem}  
\label{theoremplugin}
Suppose $X \sim p(||x - \theta||^2)$ and $Y \sim q(||y-\theta||^2), x,y, \in \mathbb{R}^d$, with $q$ absolutely continuous and strictly decreasing on $\mathbb{R}_+$. Consider estimating $q(||y - \theta||^2)$ based on $X$ under $\gamma(L(\theta, \hat{q}))$ loss,  with $\gamma$ a strictly increasing, absolutely continuous function such that $\gamma(2)<\infty$, with $\Theta$ a compact subset of 
$\mathbb{R}^d$, and with predictive density $q_{\hat{\theta},c}(y;X) = \frac{1}{c^d} \, q\left(\frac{||y - \hat{\theta}(X)||^2}{c^2}\right), y \in \mathbb{R}^d.$ Then $q_{\hat{\theta},c}(\cdot;X)$ dominates $q_{\hat{\theta},1}(\cdot;X)$ for $c \in (1,c_0)$ and some $c_0> 1$.
\end{theorem}
{\bf Proof.}   Following the proof of Theorem \ref{theoremmain} with the change of variable $y \to y - \theta$, we obtain
\begin{eqnarray*}
\nonumber 
R(\theta, q_{\hat{\theta}, c}) & = &\int\limits_{\mathbb{R}^d} \gamma \big\{ \int\limits_{\mathbb{R}^d}\, |\,q(||y - \theta||^2) - \frac{1}{c^d} q\left(\frac{||y - \hat{\theta}(x)||^2}{c^2}\right)\,|\; dy \big\} \; p(||x - \theta||^2) \; dx \\ 
 & = &  \mathbb{E}_{\theta} \, \gamma\{2\mathbb{P}_0 \big( Y \in A_1(\hat{\theta}(X)-\theta,c )\big) \, - \,  2 \, \mathbb{P}_0 \big(Y \in \mathcal{A}_2(\hat{\theta}(X)-\theta,c)\big)\} \,,
 \label{Rc22}
\end{eqnarray*}
where $A_1(\cdot,\cdot)$ and $A_2(\cdot,\cdot)$ are as defined earlier.  We also have from Lemma \ref{lemmarandomunit} that  $V \equiv \frac{Y^{\top}(\hat{\theta}(X) - \theta)}{\||Y\|\,\|\hat{\theta}(X) - \theta\|}$ is distributed independently of $(\|\hat{\theta}(X) - \theta\|, \|Y\|)$ with c.d.f. $F_V$.   Proceeding again as in Theorem \ref{theoremmain} and with the same notation, we obtain 
\begin{equation}
\label{R(c)plug-in}
R(\theta, q_{\hat{\theta},c}) = \, \mathbb{E}_{\theta} \gamma \big[ 2 \mathbb{E}_0 \left\{F_V \big( l_1^c(||\hat{\theta}(X) - \theta||,||Y||)\big) - F_V \big( l_2^c(||\hat{\theta}(X) - \theta||,||Y||)\big)\right\}\big],
\end{equation}
where the outside and inside expectations $\mathbb{E}_{\theta}$ and $\mathbb{E}_{0}$ are taken with respect to densities $p(\|x - \theta\|^2)$ and $q(\|y\|^2)$ respectively, and  $l_1^c(t_1,t_2)$  and $ l_2^c(t_1,t_2) $ are as defined in Theorem \ref{theoremmain}. Now, for fixed $\theta \in \Theta$, it follows that 
\begin{align*}
\frac{\partial}{\partial c} R(\theta, q_{\hat{\theta},c}) \Bigr|_{\substack{c \rightarrow 1^{+}}} =  -\mathbb{E}_{\theta}  \big\{\gamma' (H_1(\|\hat{\theta}(X)-\theta\|^2)) \,  \mathbb{E}_0 f_V\left(\frac{||\hat{\theta}(X)-\theta||}{2||Y||}\right)\frac{||\hat{\theta}(X)-\theta||}{||Y||} \big\} <  0\,, 
\end{align*}
which tells us that $R(\theta, q_{\hat{\theta},c}) - R(\theta, q_{\hat{\theta},1}) < 0$ for $1 < c < c_0(\theta)$.   Therefore, setting $c_0 =  \inf \{c_0(\theta), \theta \in \Theta\},$ we have that
$R(\theta, q_{\hat{\theta},c}) - R(\theta, q_{\hat{\theta},1}) < 0$ for $1 < c < c_0$ and $\theta \in \Theta$, by compactness of $\Theta$.
 \qed

The above finding and proof are unified for all unimodal $q$, all $p$, dimension $d \geq 2$, choice of $\gamma$ and of the plug-in estimator $\hat{\theta}(X)$.   A detailed illustration which expands on the determination of $c_0$ is presented in Section 4. 

\subsection{The uniform case}
\label{uniformscaleexpansion}

Notwithstanding the univariate findings of Kubokawa et al. (2017) for logconcave densities $q(|y-\theta|^2)$ for $Y$, there is, to the best of our knowledge, no compelling a priori reason why the plug-in density $q(|y-X|^2)$ should be improvable (or not improvable) by scale expansion.  In fact as mentioned in Remark \ref{remarkno}, given the optimal equivariant property of the latter density in the context of Theorem \ref{theoremuniform} and under its given conditions (i.e., $B \geq A/2$), it follows that no improvement is possible among scale expansion or scale shrinking densities $\hat{q}_c(\cdot;X)$ for such uniform models.  It is because the corresponding uniform model density is not logconcave that this example does not contradict the results of Kubokawa et al. (2017). 

It is instructive to revisit the univariate uniform model case from a more general perspective with respect to the distribution of $X$, and to study a multivariate extension which we now proceed in doing so. 

\subsubsection{Univariate case} 
In the univariate case study which follows, we set $X \sim p(|x-\theta|^2)$ and $Y \sim U(\theta-B, \theta+B)$, with known $p$ and $B>0$, and analyze the frequentist risk of predictive densities $\hat{q}_c(\cdot;X)$ taken as that of a $U(X-c, X+c)$ density.  We can assume $B=1$ without loss of generality since we can set $(X',Y')=(X/B,Y/B)$ with $X' \sim p^*(|x-\theta'|^2)$, $Y' \sim U(\theta'-1, \theta'+1)$, $\theta'=\theta/B$, and $p^*(t)\,=\, B \, p(Bt)$ for $t>0$.

\begin{theorem}
\label{theoremuniformunscaled=1}
Let $X \sim p(|x-\theta|^2)$, $Y \sim U(\theta-1, \theta+1)$, and consider the estimation of the density $\frac{1}{2} \, \mathbb{I}_{(\theta-1,\theta+1)}(y)$ of $Y$ under $L_1$ loss.  Then, among predictive densities  $\hat{q}_c(y;X) \, = \, \frac{1}{2c} \, \mathbb{I}_{(X-c,X+c)}(y)$, either one of the following conditions is sufficient for the  plug-in choice $\hat{q}_1$ to be optimal:
\begin{enumerate}
\item[ \bf{(i)}]
 $E_0(|X| \big||X|\leq 2) \leq 1$ and $p(s) \geq p(s+2)$ for all $s > 0$;

\item[ \bf{(ii)}]   $p$ is non-increasing on $(0,\infty)$.
\end{enumerate}
\end{theorem}
{\bf Proof.}   The distribution of $|X| \big| |X \leq 2$, which has density proportional to $p(s^2) \mathbb{I}_{(0,2)}(s)$, is under condition {\bf (ii)} stochastically smaller than that of a $U(0,2)$ distribution.    Therefore, it is easy to see that {\bf (ii)} implies {\bf (i)} so that we only need to establish {\bf (i)}. 

Densities $\hat{q}_c$ have constant risk so it suffices to study the risk at $\theta=0$ denoted $R(c)=R(0,\hat{q}_c)$.   For $c \in (0,1)$, the loss incurred by $\hat{q}_c$, for $\theta=0$ and as a function of $x$, becomes
\begin{eqnarray*}
\nonumber  L\big(0,\hat{q}_c(\cdot;x)\big) \,& = &\, \int_{\mathbb{R}} \big|\, \frac{1}{2c} \, \mathbb{I}_{(x-c,x+c)}(y) \, - \, \frac{1}{2} \, \mathbb{I}_{(-1,1)}(y)\, \big| \, dy \\ 
 \, & = & 2 \, \mathbb{I}_{(1+c, \infty)}(x) \, + \, (|x|-c+1) \, \mathbb{I}_{(1-c, 1+c]}(x) \, + \,  2 \, (1-c) \, \mathbb{I}_{(0, 1-c]}(x)\,.
\end{eqnarray*}
 Therefore, for $c \in (0,1)$ and denoting $f$ as the density of $|X|$, 
 \begin{eqnarray*}
 R(c) \, & = & (2-2c) \, \int_0^{1-c} f(t) \, dt \, + \, \int_{1-c}^{1+c} (t-c+1) \, f(t) \, dt \, + \, 2 \, \int_{1+c}^{\infty} f(t) \, dt\, \\
 \, & = & 2 \, + \, (1-c) \, \mathbb{P}_0 (|X| < 1-c) \, - \, (1+ c) \, \mathbb{P}_0 (|X| < 1+c) \, + \, \int_{1-c}^{1+c} t \, f(t) \, dt\,.
 \end{eqnarray*}
With derivative $R'(c) \, = \, - \left\lbrace\mathbb{P}_0(|X| < 1-c) \, + \mathbb{P}_0(|X| < 1+c) \right\rbrace < 0$ for $c \in (0,1)$, it follows that $\inf_{c>0} R(0,\hat{q}_c) \, = \, \inf_{c\geq 1} R(0,\hat{q}_c)$.    

For $c \geq 1$, we proceed in a similar fashion to obtain 
\begin{equation}
 \, 2 \, \mathbb{I}_{(c+1, \infty)}(x) \, + \, (\frac{|x|-1}{c}+1) \, \mathbb{I}_{(c-1, c+1]}(x) \, + \,  2 \,(1-1/c) \, \mathbb{I}_{(0, c-1]}(x)\,. 
\end{equation} 
   From this, we obtain for $c \geq 1$ the risk
 \begin{eqnarray*}
 R(c) \, & = & (2-\frac{2}{c}) \, \mathbb{P}_0 (|X| \leq c-1)\, + \,\frac{1}{c} \int_{c-1}^{c+1} (c-1+t) \, f(t) \, dt \, + \, 2 \, \mathbb{P}_0 (|X| > c+1)\, \\
 \, & = & 2 \, + \, (1-\frac{1}{c}) \, \mathbb{P}_0 (|X| \leq c-1) \, - \, (1+ \frac{1}{c}) \, \mathbb{P}_0 (|X| < c+1) \, + \, \frac{1}{c} \, \int_{c-1}^{c+1} t \, f(t) \, dt\,,
 \end{eqnarray*}
and its derivative $R'(c) \, = \, G(c)/c^2$ with $G(c)\,=\, \mathbb{P}_0 (|X| \leq c-1) + \mathbb{P}_0 (|X| \leq c+1) \, - \, \int_{c-1}^{c+1} t \, f(t) \, dt\,$; $f$ being the density of $|X|$.
Since $G'(c)\,=\, c \, \{f(c-1)\,-\, f(c+1)\}$, it follows under condition {\bf (i)} that $R'(c)$ changes signs at most once from $-$ to $+$ on $(0,\infty)$.  Finally, we have $\lim_{c \to 1^+}R'(c)\,=\, \mathbb{P}_0 (|X| \leq 2) \{1- E_0(|X|||X|\leq 2)\}$, and the result follows.  
\qed

\begin{example}

Theorem \ref{theoremuniformunscaled=1}'s optimality finding for the plug-in density $\hat{q}_1$ applies quite generally for non-decreasing $p$ highlighting the significance of the target uniform distribution in the formulation of the problem.   
The theorem covers cases where $X \sim U(\theta-A, \theta+A)$, with $A>0$.  Of course, Theorem \ref{theoremuniform} establishes the stronger optimality minimum risk equivariant property of $\hat{q}_1$ for $A \leq 2$, but Theorem \ref{theoremuniformunscaled=1} here applies also for $A>2$.   

Theorem  \ref{theoremuniformunscaled=1} is also applicable to samples with $X_1, \ldots, X_n$ i.i.d. $p(|x-\theta|^2)$  and plug-in densities $\hat{q}_c \sim U(\hat{\theta}-1, \hat{\theta}+1)$ where $\hat{\theta}$ is an estimator of $\theta$ based on $X_1, \ldots, X_n$.  Indeed, setting $\hat{\theta}=X$, $\hat{q}_1$ is optimal among $\hat{q}_c$'s as long as $X$ satisfies the conditions of the theorem.
Such examples are plentiful and include normal models $X_1, \ldots, X_n$ i.i.d. $N( \theta, \sigma^2)$ with $\hat{\theta}=\bar{X}$.  Another finding, which relates again to Theorem \ref{theoremuniform},  is given by  $X_1, \ldots, X_n$ i.i.d. $  U(\theta-A, \theta+A)$ with $\hat{\theta}(X) = (X_{(1)}+X_{(n)})/2$.  After checking that $|X-\theta|$ has density $\frac{2n}{(2A)^n} \, (A-2s)^{n-1} \, \mathbb{I}_{(0,A)}(s)$, one infers the optimality of the plug-in density $\hat{q}_1$ via condition {\bf (ii)} of Theorem \ref{theoremuniformunscaled=1}.

\end{example}

\subsubsection{Multivariate case}

In this section, we consider the problem of estimating the uniform density over a ball in $\mathbb{R}^d$ centered at $\theta$ based on $X \sim p(\|x - \theta\|^2)$.
We denote $\mathcal{B}_{\mu}(m)=\{t \in \mathbb{R}^d: \|t-\mu\| \leq m  \}$ as the ball of radius $m>0$ centered at $\mu \in \mathbb{R}^d$, and $\mathcal{V}_d(m)$ its volume which is equal to $m^d \pi^{d/2}/\Gamma(\frac{d}{2}+1)$.  We thus denote $Y \sim U(\mathcal{B}_{\theta}(1))$; the radius set to $1$ without loss of generality; and the associated target density is given by $q(\|y-\theta\|^2)\,=\,\frac{1}{\mathcal{V}_d(1)} \, \mathbb{I}_{\mathcal{B}_{\theta}(1)} (y)$.  The following intermediate results will be helpful to analyze the frequentist risk of scale modifications $\hat{q}_c(\cdot;X)$.
 
\begin{lemma}
\label{lemmamarginalY1}
Let $Y= (Y_1,\ldots,Y_d)^{\top} \sim U(\mathcal{B}_{0}(1))$. Then, the univariate marginal density of $Y_1$ is given by
\begin{equation}
    f_{Y_1}(t) = \frac{\Gamma\left(\frac{d}{2} +  1\right)}{\sqrt{\pi} \, \Gamma\left(\frac{d + 1}{2} \right)} \, (1 - t^2)^{\frac{d-1}{2}} \, \mathbb{I}_{(-1,1)}(t).
    \label{densityfunction_y1}
\end{equation}
\end{lemma}
{\bf Proof.}  The marginal distribution function of $Y_1$ is given by
\begin{align}
   \nonumber  F_{Y_1}(t) &  = \frac{1}{\mathcal{V}_d(1)} \int\limits_{-1}^{t} \mathcal{V}_{d-1}(\sqrt{1 - s^2}) \;ds \\
    & = \int\limits_{-1}^{t} \frac{\Gamma(\frac{d}{2} + 1 )}{\sqrt{\pi}\Gamma(\frac{d+1}{2})} (1 - s^2)^{\frac{d-1}{2}}\; ds \text{, for } t \in (-1,1),
\end{align} 
which is indeed the c.d.f. of density (\ref{densityfunction_y1}).  \qed

\begin{lemma}
\label{lemmaintersectionballs}
For all $d \geq 2$, $x \in \mathbb{R}^d-\{0\}$, $c>0$, the intersection of the balls  $\mathcal{B}_x(c)$ and $\mathcal{B}_0(1)$ has volume
\begin{equation}
\label{intersectionballs}
Vol \big( \mathcal{B}_x(c) \cap \mathcal{B}_0(1) \big) \, = \, \mathcal{V}_d(1) \, \left\lbrace1- F_{Y_1} \big(\frac{\|x\|}{2}  + \frac{ 1- c^2 }{2\|x\|}  \big) \, + \, c^d F_{Y_1} \big(-\frac{\|x\|}{2c}  + \frac{ 1-  c^2 }{2 \, c \, \|x\|}  \big)\right\rbrace\,,
\end{equation}
where $F_{Y_1}$ is the c.d.f. given in Lemma \ref{lemmamarginalY1}.
\end{lemma}
{\bf Proof.}  Observe that the given ratio of volumes depends on $x$ only through its norm $\|x\|$ so that, without loss of generality, we can set $x_0 = (\|x\|,0,\ldots,0)^{\top}$. The points $z = (z_1,z_2^{\top})^{\top} \in \mathbb{R}^d$ with $z_2 \in \mathbb{R}^{d-1}$ lying on the intersection of $\mathcal{B}_{x}(c)$ and $\mathcal{B}_{0}(1)$  satisfy the equations $\|z\|^2 = 1$ and $\|z-x_0\|^2 = c^2$, i.e., $z_1 = \frac{||x||}{2}  - \frac{ c^2 - 1}{2||x||}$ and $\|z_2\|^2 = 1 - z_1^2.$ We thus obtain
\begin{align*}
     \frac{Vol\big(\mathcal{B}_{x}(c) \cap  \mathcal{B}_{0}(1)\big)}{\mathcal{V}_d(1)} &   = \frac{Vol(\left\{z \in \mathbb{R}^d : \|z\|^2 \leq 1 \;,\; z_1 \geq \frac{\|x\|}{2}
      - \frac{ c^2 - 1}{2\|x\|} \right\})}{\mathcal{V}_d(1)} \\ & + \frac{Vol(\left\{z \in \mathbb{R}^d : \|z - x\|^2 \leq c^2 \;,\; z_1 \leq \frac{\|x\|}{2}  - \frac{ c^2 - 1}{2\|x\|} \right\})}{\mathcal{V}_d(1)} \\
      & = \mathbb{P}\left(Y_1 \geq  \frac{\|x\|}{2}
      - \frac{ c^2 - 1}{2\|x\|} \right) \, + \, c^d \, \mathbb{P}\left( cY_1 + \|x\| \leq \frac{\|x\|}{2}
      - \frac{ c^2 - 1}{2\|x\|} \right),
\end{align*}
which yields the result. 
\qed

The next result provides the risk of $\hat{q}_c$.
\begin{theorem}
\label{theoremriskuniformulamultivariate}
Consider $X \sim p(\|x - \theta\|^2)$, $Y \sim U(\mathcal{B}_{\theta}(1))$   independent, and consider estimating the density of $Y$ under $L_1$ loss.    Then, the risk of the predictive density  $\hat{q}_{c}(y;X)  = \ \frac{1}{c^d\mathcal{V}_d(1)} \mathbb{I}_{\mathcal{B}_{X}(c)}(y)$  is given by
\begin{align}
\nonumber R(c ) & = 2 \left[ \mathbb{E}^{\|X\|} \left\{ F_{Y_1}\left( \frac{||X||}{2}  - \frac{ c^2 - 1}{2||X||}\right)  + c^d F_{Y_1}\left(  \frac{||X|| }{2c} +  \frac{c^2 - 1}{2c||X||}\right) \right\} - c^d\right]   & \text{ for } c \leq 1, \\
\nonumber & = \frac{2}{c^d} \left[ \mathbb{E}^{\|X\|} \left\{ F_{Y_1}\left( \frac{||X||}{2}  - \frac{ c^2 - 1}{2||X||}\right)  + c^d F_{Y_1}\left(  \frac{||X|| }{2c} +  \frac{c^2 - 1}{2c||X||}\right) \right\} - 1\right] & \text{ for } c > 1 .
\end{align}

\end{theorem}
{\bf Proof.} 
It is interesting here to exploit the relationship between $L_1$ loss (say $\rho(q_{\theta},\hat{q}_{c})$) and the overlap coefficient.   As in (\ref{formularisk}), the risk of $\hat{q}_{c}(y;X)$ is constant in $\theta$ and given by

\begin{equation}
\label{formularisktheorem3.5}
R(c) \, = \, \int\limits_{\mathbb{R}^d} \rho\big( q_0, \hat{q}_c(\cdot;x) \big) \; p(||x||^2) \; dx\,,
\end{equation}
with $q_0(y) \, = \, \frac{1}{\mathcal{V}_d(1)} \, \mathbb{I}_{\mathcal{B}_0(1)}(y)$.  From (\ref{overlap}), we obtain
\begin{eqnarray*}
\rho\big(q_0, \hat{q}_c(\cdot;x)\big) & \, = \, & 2 \, - 2 \, \int_{\mathcal{B}_{x}(c) \cap  \mathcal{B}_{0}(1))} \; \min\big(q_0(y), \hat{q}_c(y;x)\big) \, dy \\
\, & \, = \, &  2 \, - \, 2 \, \frac{Vol\big(\mathcal{B}_{x}(c) \cap  \mathcal{B}_{0}(1)\big)}{\max(1, c^{d})\, \mathcal{V}_d(1)}\,,
\end{eqnarray*}
and the result follows from (\ref{formularisktheorem3.5}) and Lemma \ref{lemmaintersectionballs}.  \qed

\begin{remark} 
\label{expressionR(1)}
Observe that the value $R(1)$, i.e., the $L_1$ risk of the plug-in density $\hat{q}_1(\cdot;X)$, reduces to 
\begin{equation}
\label{97}
4 \;\mathbb{E}^{||X||} F_{Y_1}\left(\frac{||X||}{2}\right) - 2\,.
\end{equation}
This matches the expression given by Kubokawa et al. (2017), established under the conditions of Theorem \ref{theoremmain}, with $F_{Y_1}$ the c.d.f. of $Y_1$ for  $Y=(Y_1, \ldots, Y_d)^{\top} \sim q(\|y\|^2)$ and $q$ strictly decreasing.  Interestingly, a direct reading of $R(1)$ from (\ref{R(c)gamma}) with identity $\gamma$ yields the expression
\begin{equation}
\label{96}
R(1) \, = \, 4 \mathbb{E}^{\|X\|} \big\{ \mathbb{E}^{\|Y\|}  F_V \big(\frac{\|X\|}{2 \|Y\|} \big)\big\} \, - \, 2\,,
\end{equation}
with $F_V$ the c.d.f. associated to (\ref{densityV}). But, expressions (\ref{97}) and (\ref{96}) are seen to match with the independence of $X$ and $Y$, and since for all $a>0$:
\begin{equation}
\nonumber
F_{Y_1}(a) \, = \, \mathbb{E}^{\|Y\|} \, \mathbb{P} \left( \frac{Y_1}{\|Y\|} \leq \frac{a}{\|Y\|}  \, \big| \, \|Y\|\right) \, = \, 
\mathbb{E}^{\|Y\|} \, F_V \big( \frac{a}{\|Y\|} \big)\, ,
 \end{equation}
given that $\frac{Y_1}{\|Y\|}$ and $\|Y\|$ are independently distributed with $\frac{Y_1}{\|Y\|} =^d V$ (Lemma \ref{lemmarandomunit}).
\end{remark}

\begin{figure}
  \centering
  \includegraphics[width=0.9\linewidth]{"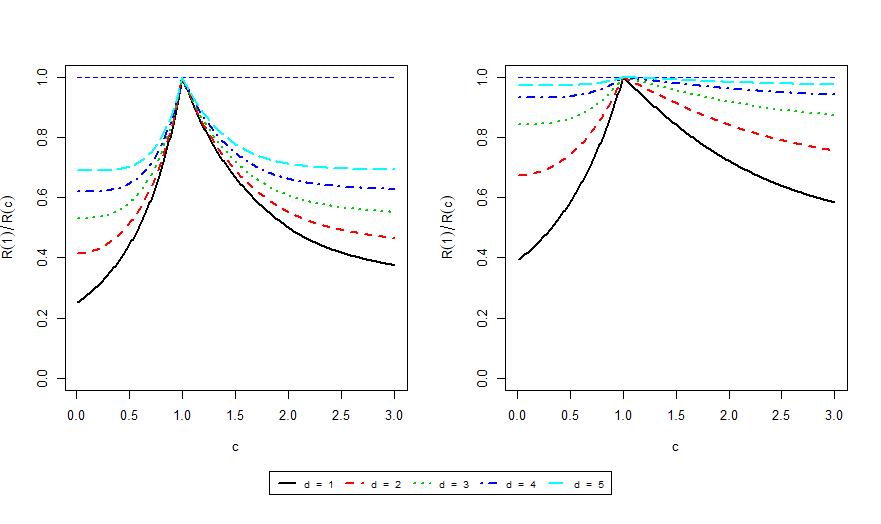"}
\caption{Risk ratios $\frac{R(1)}{R(c)}$ as a function of $c$ for $d = 1,2,3,4,5$, $Y \sim U(\mathcal{B}_{\theta}(1))$,  with $X \sim U(\mathcal{B}_d(\theta,1))$ on the left and $X \sim N_d(\theta, I_d)$ on the right.}
\label{uniformden}
\end{figure}

\normalsize
We proceed with further analysis of $R(c)$, including numerical evaluations presented in {\bf Figure  \ref{uniformden}}, which exhibit the plausible optimality of $\hat{q}_1$ for $d = 2,3,4,5$ among $\hat{q}_c$ for $X \sim U(\mathcal{B}_{\theta}(1))$ and  $X \sim N_d(\theta,I_d)$, the result having been established earlier in Theorem 3.4 for $d = 1.$ Exact analysis is difficult to achieve, but we do conclude this section with a definite answer for $d=3$ and $X \sim {U}(\mathcal{B}_{\theta}(1))$.

\begin{theorem}
\label{theoremuniformunscaled=3}
Suppose $d =3, $ and  $Y \sim U(\mathcal{B}_{\theta}(1))$ and $X \sim U(\mathcal{B}_{\theta}(1))$, independent and consider the estimation of the density $\frac{3}{4\pi} \, \mathbb{I}_{\mathcal{B}_{\theta}(1)}(y)$ of $Y$ under $L_1$ loss.  Then, among predictive densities  $\hat{q}_c(y;X) \, = \, \frac{3}{4c^3\pi} \, \mathbb{I}_{\mathcal{B}_{X}(c)}(y)$,  the  plug-in choice $\hat{q}_1$ is optimal.
\end{theorem}
{\bf Proof.}  For $d=3$, we have $F_{Y_1}(t)\,=\, (3t-t^3+2)/4$.  For $X \sim U(\mathcal{B}_{\theta}(1))$, the density function of $R =\|X\|$ is given by $h(r) = 3r^2\, \mathbb{I}_{(0,1)}(r)$.   With these expressions, an evaluation of the risk in Theorem \ref{theoremriskuniformulamultivariate} gives
\begin{equation}
    R(c) = \frac{c(18-c^2)}{16} \mathbb{I}_{(1,2]}(c) \, + \, 2\left(1 - \frac{1}{c^3}\right) \mathbb{I}_{(2,\infty)}(c) + \frac{-c^6 + 18c^4 - 32c^3 + 32}{16} \, \mathbb{I}_{(0,1)}(c).
\end{equation}
It is then simple to verify that $R(c)$ is decreasing on $(0,1)$ and increasing on $(1,\infty)$, which establishes the result (with $\inf_c R(c) \, = \, R(1) \,=\, 17/16$).

\section{Examples}

In this section, we illustrate the dominance findings in terms of implementation and with frequentist risk comparisons.  The dominance results are wide ranging with respect to the model specifications $p$ and $q$, the dimension $d$, the choice of loss $\gamma$, and the choice of the plug-in density for the compact parameter space case, but we focus nevertheless on two specific situations involving normal models.  In the first case, we compare in the context of Section 3.1 the risks of the optimal $\hat{q}_{c^*}$ with that of the plug-in density $\hat{q}_1$.  In the second case, we illustrate the findings of Section 3.2 for a mean parameter restricted to a ball.

\begin{example}
Consider $X \sim N_d(\theta,\sigma_X^2 \, I_d)$ and $Y \sim N_d(\theta,\sigma_Y^2 I_d)$, $d \geq 2$, and $L_1$ loss.  Theorem \ref{theoremmain} tells us that there exists, among scale modifications $\hat{q}_c$, an optimal $\hat{q}_{c^*}$ with $c^*>1$.   The optimal value of $c^*$ can be evaluated numerically by evaluating $R(c)$, which in turn is given by (\ref{R(c)gamma})
with identity $\gamma$, $\ell_1^c(t_1,t_2)  = \frac{t_1^2 - (c^2 - 1)t_2^2 + 2dc^2\sigma^2_Y\text{log}c}{2t_1t_2}$ and $\ell_2^c(t_1,t_2) = \frac{-t_1^2 - (c^2 - 1)t_2^2 + 2d\sigma^2_Y 
\log c}{2ct_1t_2} $. The optimal value $c^*$ depends on the dimension $d$ and can be shown to depend on $(\sigma^2_{X}, \sigma^2_Y)$ only through the ratio $r= \sigma^2_Y/\sigma^2_X.$  For fixed $r$, $c^*$ turns out to be remarkably stable as a function of $d$, increasing slightly towards a limiting value at $d \to \infty$ as exhibited in {\bf Figure \ref{OptimalcRiskratios}}.    It would be interesting to identify further analytical properties relative to $c^*$, in particular as $d \to \infty$; but we have been unable to do so.  Such properties will depend on the underlying spherically symmetric model.  Further numerical evidence which is not portrayed here suggests quite different behaviour, with for instance $c^*$ diverging for a multivariate Cauchy model for $X$ and $Y$.     

In terms of risk, the plug-in density achieves the constant risk
\begin{equation}
\label{equationR(1)}
 R(1) \, = \, 4 \mathbb{E} \Phi \big(\frac{\|Z\|}{2 \sqrt{r}} \big) \, - \, 2  \, , \hbox{ with } Z \sim N_d(0,I_d)\,,
\end{equation}
 as noted in Remark \ref{expressionR(1)}.  As a complementary note, we provide in the Appendix a nice expansion for $R(1)$ and it namely simplifies to $R(1)\,=\, \frac{2}{\sqrt{1+4r}}$ for $d=2$.
The risk $R(1)$ in (\ref{equationR(1)}) clearly decreases in $r$, and increases in $d$ given the stochastic increasing ordering of $\|Z\|.$ The increasingness in terms of the ratio of the variances $r$ is associated with a relative better concentration of $X$ about $\theta$, and consequently more precise estimates $\hat{q}_1(\cdot;X)$ for estimating $q_{\theta}$.   The frequentist risk of $\hat{q}_{c^*}$ will also decrease as a function of $r$, and the relative merits of these two densities are represented in {\bf Figure \ref{OptimalcRiskratios}} with graphs of the ratio $\frac{R(c^*)}{R(1)}$ for varying $d$ and for $r=1/2, 1, 2$.    The gains are modest, they are attenuated as $d$ increases, but do not exhibit a straightforward ordering with changes in $r$.

\begin{figure}
  \centering
  \includegraphics[width=0.9\linewidth]{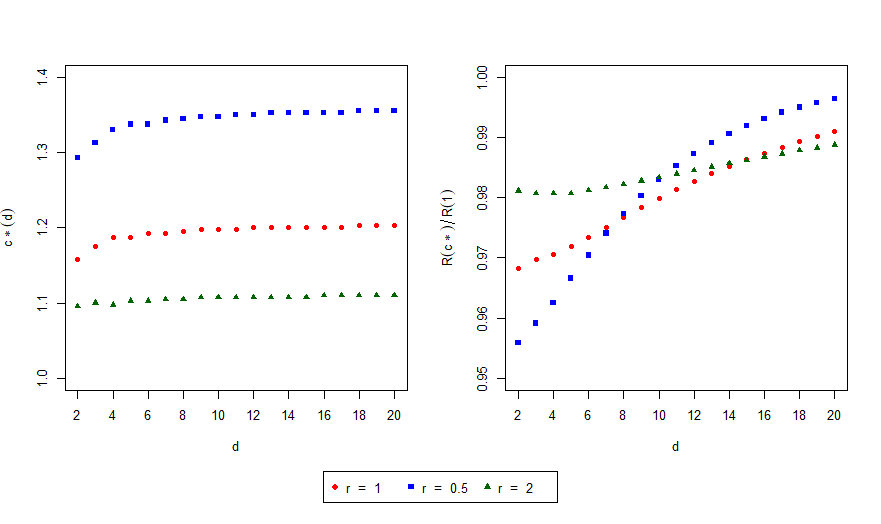}
\caption{Optimal $c^*$ values and risk ratios $\frac{R(c*)}{R(1)}$ as functions of $d$ for normally distributed $X$ and $Y$. }
\label{OptimalcRiskratios}
\end{figure}

\end{example}

\begin{example}
\label{exampleball}
We illustrate our previous finding for the case where $\theta$ is constrained to the  ball of radius $m$ around the origin (i.e., $\theta \in \mathcal{B}_0(m))$.
Simply stated, for any non-degenerate estimator $\hat{\theta}(X)$ of $\theta$, $d \geq 2$, and unimodal $q$, Theorem \ref{theoremplugin} says that is inevitable that there will be scale expansion variants $q_{\hat{\theta},c} $ that dominate the plug-in density $q_{\hat{\theta},1} $ under $L_1$ frequentist risk for $\theta \in B_m$.   We refer to Marchand \& Strawderman (2004) and references therein for point estimation aspects under such parametric restrictions. 

For the purpose of illustration, consider $X \sim N_{d}(\theta, I_d)$, $Y \sim N_{d}(\theta, I_d)$ where $ \lambda = \|\theta\| \leq m.$  Given the parametric constraint, the maximum likelihood estimator $\hat{\theta}_{mle}(X) = \min\{m , \|X\|\} \, X/\|X\|$ is an appealing choice to estimate $\theta$, and we thus consider predictive densities  $q_{\hat{\theta}_{mle},c}(y;X) \, = \, 
\frac{1}{c^d} \, q \big( \frac{\|y-\hat{\theta}_{mle}(X)\|^2}{c^2}  \big)$; $y \in \mathbb{R}^d$; with $c=1$ the plug-in mle density, and $c>1$ yielding scale expansion variants.   In the spirit of Theorem \ref{theoremplugin}'s strategy of proof, one can approach numerically: (i) for fixed $\lambda \in [0,m]$, the optimal value $c^*_1(\lambda)$ minimizing the risk for $\|\theta\|=\lambda$ of $q_{\hat{\theta}_{mle},c}$, and then (ii) the value of $c_1=\inf_{\lambda \in [0,m]} \{c^*_1(\lambda)\}$ with dominance on the range $c \in (1,c_1]$ a consequence of Theorem \ref{theoremplugin}. The obtained $c_1$ will be a lower bound for Theorem 3.3's $c_0$ value. The parametric reduction arises as the $L_1$ risk depends on $\theta$ only through $\lambda$, and this is expanded upon in the Appendix.  

{\bf Figure \ref{pluginMLE}} presents for $d=3$ and $m=1$ the frequentist risks, as functions of $\lambda$, of the predictive densities
$q_{\hat{\theta}_{mle,1.05}} \sim N_3(\hat{\theta}_{mle}(X), (1.05)^2 I_3)$, $q_{\hat{\theta}_{mle,1}} \sim N_3(\hat{\theta}_{mle}(X), I_3)$, along with $N_3(X, I_3)$ and $\hat{q}_{c^*} \sim N_3(X, (c^*)^2 I_3)$ densities.   The scale expansion level $c=1.05$ was chosen in accordance with the above strategy with $c_1 \approx 1.052$. 
 The theoretical improvement for $0 \leq \lambda \leq 1$ of $\hat{q}_{mle,c_1}$ over $\hat{q}_{mle,1}$ is very slight, but still present.  Based on 
numerical evidence, the dominance persists for the misspecification $\lambda >1$.    This is illustrated in Figure 3 for values $\lambda \in (1,2.2)$, but also is inferred
by calculations of $c^*_1(\lambda)$, with $c^*_1(\lambda)>1.05$ for all $\lambda >1$. Both of these choices are largely better that the others, and even offer improvement for some level of misspecification, with $R(\theta,q_{\hat{\theta}_{mle,1.05}}) \leq R(\theta, \hat{q}_{c^*}) \hbox { iff } \|\theta\| \leq k $ with $k \approx 2.1$.   Both the $N_3(X, (c^*)^2I_3)$,
and $N_3(X, I_3)$ densities ignore the parametric restriction, and the latter ignores as well scale expansion improvement, which is optimal here for $c^* \approx 1.175$ (see {\bf Figure \ref{OptimalcRiskratios}}).    

We point out that the dominance of $q_{\hat{\theta}_{mle},1}$, and therefore of $q_{\hat{\theta}_{mle},c_1}$, over $\hat{q}_0$ is theoretically justified since the general $L_1$ risk comparison of plug-in densities relates directly to the point estimation risk comparison under loss 
\begin{equation}
\label{lossdual}
4 \, \Phi \big( \frac{\|\hat{\theta}-\theta\|}{2} \big) \, - \, 2 \,
\end{equation}
(i.e., Corollary 2.1 of Kubokawa et al., 2017).
Since $\|\hat{\theta}_{mle}(X)-\theta\|^2$ is stochastically smaller than $\|X-\theta\|^2$ for all $\theta \in \Theta(m)$, it follows that $\hat{\theta}_{mle}(X)$ dominates $\theta_0(X)$ under loss (\ref{lossdual}), and therefore that $q_{\hat{\theta}_{mle,1}}$ dominates $\hat{q}_0$ under $L_1$ loss.

\begin{figure}
  \centering
  \includegraphics[width=0.75\linewidth]{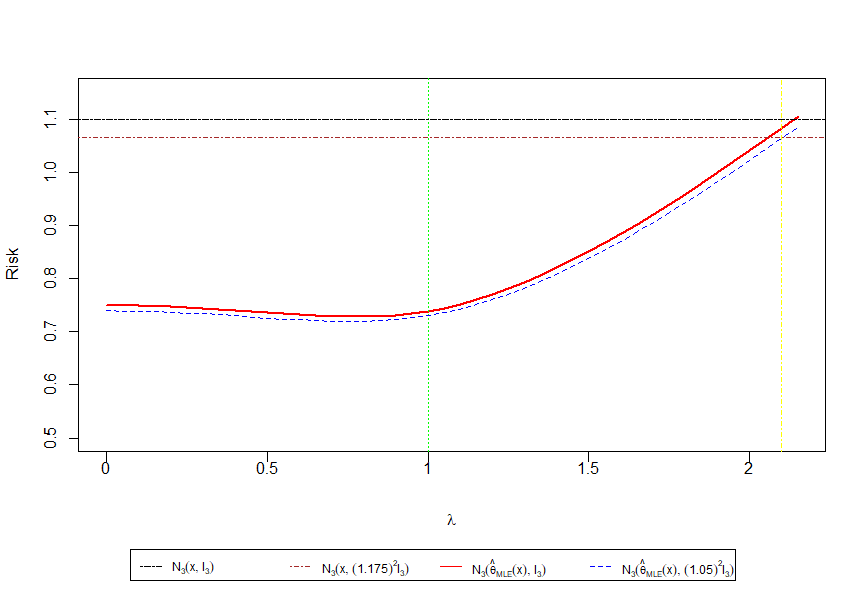}
  \caption{Frequentist risks as function of $\lambda = \|\theta\|$ for various predictive densities of $Y | \theta \sim N_3(\theta, I_3)$ based on $X | \theta \sim N_3(\theta, I_3)$ with $\|\theta\| \leq m = 1$. }
\label{pluginMLE}
\end{figure}

\end{example}

\section{Concluding remarks}

With this paper, we have addressed technical challenges present in the $L_1$ risk analysis of predictive densities, focussing on spherically symmetric models.   We have provided improved scale expansion variants of the plug-in density $\hat{q}_0(y;X) \, = \, q(\|y-X\|^2)$, $y \in \mathbb{R}^d$, for estimating the density of $Y \sim q(\|y-\theta\|^2)$ based on $X \sim  p(\|x-\theta\|^2)$.   We have shown that such improvement is inevitable quite generally with respect to the specifications of $p$, unimodal $q$, and $d>1$; and also elaborated on the necessity of the unimodality assumption.  The findings complement those of Kubokawa et al. (2017), namely in obtaining dominating predictive densities of $\hat{q}_0(\cdot;X)$ for $d=2,3$.  Furthermore, we have obtained novel extensions to losses which are increasing functions of $L_1$ loss, and also to other plug-in densities quite generally when the parameter space is compact.

Despite the natural appeal of $L_1$ distance, including its relationships to total variation distance and the overlap coefficient, there have been quite few previous results in the literature and the challenges met here could well provide avenues for future work, namely for spherically symmetric models with unknown scale with or without a residual vector (e.g. Kato, 2009; Boisbunon \& Maruyama, 2014; Fourdrinier et al., 2019), and for non-symmetric multivariate models such as those described by skewed multivariate normal densities.

\section*{Acknowledgements}
\'{E}ric Marchand's research is supported in part by the Natural Sciences and Engineering Research Council of Canada.   Pankaj Bhagwat is grateful to the ISM (Institut des sciences math\'ematiques) for financial support.  We are grateful to Bill Strawderman for fruitful discussions in particular on Bayesian perspectives for the uniform model studied in Section 2.2.   We are thankful to Benjamin Heuclin for useful and insightful preliminary numerical evaluations.

\section*{Appendix}
{\bf Lemma 5.5}

For $X \sim N_d(\theta, \sigma^2_X I_d)$ and $Y \sim N_d(\theta, \sigma^2_Y I_d)$, an explicit expression for the $L_1$ risk of density $N_d(X, \sigma^2_Y I_d)$ is given by 
\begin{equation} 
\label{L1lossunbiased}
 R(1) \, = \, \frac{\Gamma\big((d+1)/2\big)}{\Gamma(d/2)} \sqrt{\frac{4}{r \pi}} \, _2F_1(\frac{d+1}{2}, \frac{1}{2}, \frac{3}{2}; -\frac{1}{4r})\,,
\end{equation}
where $r=\frac{\sigma_Y^2}{\sigma_X^2}$ and $_2F_1$ is Gauss' hypergeometric function.  For the particular case of $d=2$, this reduces to $\frac{2}{\sqrt{1+4r}}$.

\noindent  {\bf Proof.} The $d=2$ case follows as $\,_1F_0(\frac{1}{2};-;t)\,=\,(1-t)^{-1/2}$ for $t<1$. To establish (\ref{L1lossunbiased}), we start with expression
(\ref{equationR(1)}) which we write as  $R(1) \, = \, 4 \mathbb{E} \, \Phi (\sqrt{U}) \, - \, 2\,$ with $U \sim Ga(\frac{d}{2}, \frac{1}{2r})$.   An expansion yields for $u \in \mathbb{R}$
\begin{eqnarray*}
\Phi(\sqrt{u}) \, & = & \, \frac{1}{2} \, + \, (2\pi)^{-1/2} \, \int_0^{\sqrt{u}}  \sum_{j \geq 0} 
\frac{(-1)^j  \, (\frac{t^2}{2})^j} {j!} \, dt \\
\, & = & \frac{1}{2} \, + \, (2\pi)^{-1/2} \, \sum_{j \geq 0} \frac{(-1)^j \, u^{j+1/2}}{j! \, 2^j \, (2j+1)}\, \\
\, & = & \frac{1}{2} \, + \, (2 \pi)^{-1/2} \, \sum_{j \geq 0} \frac{(-1)^j \, u^{j+1/2} (1/2)_j}{j! \, 2^j \, (3/2)_j}\,.
\end{eqnarray*}
From this, expression (\ref{L1lossunbiased}) follows by taking expectation, extracting the moments $\mathbb{E}(U^{j+1/2}) \, = \, \frac{\Gamma(\frac{d+1}{2}+j)}{\Gamma(\frac{d}{2})} (\frac{1}{2r})^{j+1/2}$, and collecting terms.  \qed
\\

{\bf On the invariance property in Example \ref{exampleball}}

We point out that, for $X \sim p(\|x-\theta\|^2)$, $Y \sim q(\|y-\theta)\|^2)$, an equivariant estimator $\hat{\theta}(X)$ with respect to orthogonal transformations, which satisfies the relationship $\hat{\theta}(Hx) \, = \, H \hat{\theta}(x)$ for $x \in \mathbb{R}^d$ and orthogonal $H$ (or equivalently which is of the form $\hat{\theta}(X) \, = \, g(\|X\|) \, X$ for some $g$), that the $L_1$ frequentist risk of a predictive density $f(\|y-\hat{\theta}(X)\|^2), y \in \mathbb{R}^d$, depends on $\theta$ only through its norm $\lambda=\|\theta\|$.  This is the case namely for plug-in densities ($f = q$) and their scale expansion variants with $f(t) = \frac{1}{c^d} \, q(\frac{t}{c^2})$.  The property holds as, for any orthogonal matrix $H$ and $\theta \in \mathbb{R}^d$, and since  $\hat{\theta}(X)$ is equivariant:
\begin{eqnarray*}
R(H\theta, f) \, & = & \, \int_{\mathbb{R}^{2d}}  p(\|x-H\theta\|^2) 
\left\lbrace \big| f(\|y-\hat{\theta}(x)\|^2) \, - \,  q(\|y-H\theta\|^2)   \big| \right\rbrace  \, dy  \, dx \\
\, & = & \, \int_{\mathbb{R}^{2d}}  p(\|Hx-H\theta\|^2) 
\left\lbrace \big| f(\|Hy-\hat{\theta}(Hx)\|^2) \, - \,  q(\|Hy-H\theta\|^2)   \big| \right\rbrace  \, dy  \, dx \\
\, & = & R(\theta, f)\,,
\end{eqnarray*}
by the transformation $(x,y) \to (H^{\top}x, H^{\top}y)$, since $(p(\|Ht\|^2),q(\|Ht\|^2)) = ((p(\|t\|^2),q(\|t\|^2))$.  \footnote{This is more generally true for the class of divergences of the form $\rho(f,g) \, = \, \int_{\mathbb{R}^d} h(\frac{f(t)}{g(t)}) g(t) \, dt$, which include Kullback-Leibler, reverse Kullback-Leibler and $\alpha$-divergence. }
\small
\section*{References}
Aitchison, J. (1975). Goodness of prediction fit. {\it Biometrika}, {\bf 62}, 547-554.

\medskip\noindent Boisbunon, A. \& Maruyama, Y.  (2014).  Inadmissibility of the best equivariant density in the unknown variance case.   {\it Biometrika}, {\bf 101}, 733-740.

\medskip\noindent Brown, L.D., George, E.I., \&  Xu, X. (2008).
Admissible predictive density estimation.
{\it Annals of Statistics}, {\bf 36}, 1156-1170.

\medskip\noindent  DasGupta, A. \& Lahiri, S.N. (2012).  Density estimation in high and ultra dimensions, regularization, and the $L1$ asymptotics.
{\it Contemporary Developments in Bayesian analysis and Statistical Decision Theory: 
A Festschrift for William E. Strawderman}, IMS Collections, {\bf 8}, 1-23.

\medskip\noindent  Devroye, L.  \& Györfi, L. (1985).  {\it Nonparametric density estimation. 
The L1 view},  Wiley, New York.

\medskip\noindent  Eaton, M.L. (1989).  Group Invariance Applications in Statistics, in: Regional Conference Series in Probability and Statistics, vol. 1, Institute of Mathematical Statistics and the American Statistical Association.

\medskip\noindent Fourdrinier, D., Strawderman, W.E. \& Wells, M. T. (2018).  {\it Shrinkage estimation}.  Springer.

\medskip\noindent Fourdrinier, D., Marchand, \'{E}., Righi, A. and Strawderman, W.E. (2011).
On improved predictive density estimation with parametric constraints.
{\it Electronic Journal of Statistics}, {\bf 5}, 172-191.

\medskip\noindent Fourdrinier, D., Strawderman, W.E. \& Wells, M. T. (2018).  {\it Shrinkage estimation}.  Springer.

\medskip\noindent
Fourdrinier, D.,  Marchand, \'E. \& Strawderman,
W.E. (2019). On efficient prediction and predictive density estimation for spherically symmetric models. {\it Journal of Multivariate Analysis}, {\bf 173}, 18-25. 

\medskip\noindent  George, E., Marchand, \'{E}., Mukherjee, G. \& Paul, D. (2019). New and evolving roles of shrinkage in large-scale prediction and inference. BIRS Workshop Report.

\medskip\noindent George, E. I., Liang, F. \& Xu, X. (2006). 
Improved minimax predictive densities under Kullback-Leibler loss. 
{\it Annals of Statistics}, {\bf 34}, 78-91.

\medskip\noindent  Kariya, T. \& Eaton M.L. (1977).  Robust tests for spherical symmetry.  {\it The Annals of Statistics}, {\bf 5}, 206--215.

\medskip\noindent Kato, K. (2009).  Improved prediction for a multivariate normal distribution with unknown mean and variance.  {\it Annals of the Institute of Statistical Mathematics}, {\bf 61}, 531-542.

\medskip\noindent Kiefer, J. (1957).  Invariance, minimax sequential
estimation, and continuous time processes. {\it Annals of
Mathematical Statistics}, {\bf 28}, 573-601.

\medskip\noindent
Komaki, F. (2001). 
A shrinkage predictive distribution for multivariate normal observables. 
{\it Biometrika}, {\bf 88}, 859-864.

\medskip\noindent Kubokawa, T., Marchand, \'{E}. \& Strawderman, W.E. (2017).
On predictive density estimation for location families under integrated absolute value loss.  {\it Bernoulli}, {\bf 23}, 3197-3212.

\medskip\noindent Kubokawa, T., Marchand, \'{E}. \& Strawderman, W.E. (2015).
On predictive density estimation for location families under integrated $L_2$ loss.   {\it Journal of Multivariate Analysis},  {\bf 142}, 57-74.

\medskip\noindent  LMoudden, A. \& Marchand, \'{E}.  (2019).
On predictive density estimation under $\alpha$-divergence loss.  {\it Mathematical Methods of Statistics},  {\bf 28}, 127--143.

\medskip\noindent Marchand, \'E., and Strawderman, W.E. (2004). {\it
Estimation in restricted parameter spaces: A review.} A
Festschrift for Herman Rubin, IMS Lecture Notes-Monograph Series
{\bf 45}, Institute of Mathematical Statistics, Hayward, CA, pp.
21-44.

\medskip\noindent Nogales, A. G. (2021).   On Bayesian estimation of densities and sampling distributions: The posterior predictive distribution as the Bayes estimator.   {\it Statistica Neerlandica}, {\bf 76}, 236--250. 

\medskip\noindent  Weitzman, M.S. (1970).   Measures of overlap of income distributions of white and negro families in the United States.    Technical Report 22, US Department of Commerce.

\end{document}